\documentclass{mcom-l}

\usepackage{amsmath}
\usepackage{amssymb}
\usepackage{graphics}
\usepackage{latexsym}

\newcommand{\pFq} {^{\phantom{(}}_pF_q^{\mbox{\tiny$(\alpha)$}}}

\newcommand{\mpFq} {^m_pF_q^{\mbox{\tiny$(\alpha)$}}}

\newcommand{\al}{{\mbox{\tiny$(\alpha)$}}}

\newcommand{\tr}{{\rm tr}}

\newcommand{\bmat}{\left[ \begin{array}}
\newcommand{\emat}{\end{array} \right]}
\newcommand{\ignore}[1]{}

\newtheorem{theorem}{Theorem}[section]
\newtheorem{lemma}[theorem]{Lemma}

\theoremstyle{definition}

\newtheorem{algorithm}[theorem]{Algorithm}

\theoremstyle{remark}

\numberwithin{equation}{section}

\begin{document}

\title[Evaluation of the Hypergeometric
Function of a Matrix Argument]
{The Efficient Evaluation of the Hypergeometric
Function of a Matrix Argument}

\author{Plamen Koev}
\address{Department of Mathematics, Massachusetts Institute of
Technology, Cambridge, Massachusetts 02139}
\email{plamen@math.mit.edu}
\thanks{This work was supported in part by NSF Grant
No.\ DMS-0314286.}

\author{Alan Edelman}
\address{Department of Mathematics, Massachusetts Institute of
Technology, Cambridge, Massachusetts 02139}
\email{edelman@math.mit.edu}


\date{\today}



\begin{abstract}

We present new algorithms that efficiently approximate the hypergeometric 
function of a matrix argument through its expansion as a series of Jack 
functions.  Our algorithms exploit the combinatorial properties of the 
Jack function, and have complexity that is only linear in the size of the 
matrix.

\end{abstract}

\maketitle

\section{Introduction}

The hypergeometric function of a matrix argument has a wide area of 
applications in multivariate statistical analysis~\cite{muirhead}, random 
matrix theory~\cite{forrester}, wireless communications 
\cite{gaosmithclark98,kangalouni04}, etc.  Except in a few special cases, 
it can be expressed only as a series of multivariate homogeneous 
polynomials, called Jack functions. This series often converges very 
slowly~\cite[p.\ 390]{muirhead78, muirhead}, and the cost of the 
straightforward evaluation of a single Jack function is 
exponential~\cite{demmelkoevschur}. The hypergeometric function of a 
matrix argument has thus acquired a reputation of being notoriously 
difficult to approximate even in the simplest cases 
\cite{butlerwood02,grs00}.

In this paper we present new algorithms for approximating the value of the 
hypergeometric function of a matrix argument. We exploit recursive 
combinatorial relationships between the Jack functions, which allow us to 
only {\em update} the value of a Jack function from other Jack functions 
computed earlier in the series. The savings in computational time are 
enormous; the resulting algorithm has complexity that is only {\em linear} 
in the size of the matrix argument. In the special case when the matrix 
argument is a multiple of the identity, the evaluation becomes even 
faster.

We have made a MATLAB~\cite{matlab} implementation of our algorithms 
available~\cite{koevwebpage}. This implementation is very efficient (see 
performance results in Section \ref{sec_numexp}), and has lead to 
new results~\cite{absiledelmankoev,edelmansutton04}.

The hypergeometric function of a matrix 
argument is defined as follows. Let $p\ge 0$ and $q\ge 0$ be integers, 
and let  
$X$ be an $n\times n$ complex symmetric matrix with eigenvalues 
$x_1,x_2,\ldots,x_n$. Then
%
%
\begin{equation}
\pFq(a_1,\ldots,a_p;
b_1,\ldots,b_q;X) \\
\equiv\sum_{k=0}^\infty\sum_{\kappa\vdash k}
\frac{(a_1)^{\al}_\kappa\cdots(a_p)_\kappa^{\al}}
{k!(b_1)_\kappa^{\al}\cdots(b_q)_\kappa^{\al}} \cdot
C_\kappa^{\al}(X),
\label{hg1}
\end{equation}
where $\alpha>0$ 
is a parameter; $\kappa\vdash k$ means $\kappa=(\kappa_1,\kappa_2,\ldots)$ 
is a partition of $k$ (i.e., $\kappa_1\ge\kappa_2\ge\ldots\ge 0$ are 
integers such that $|\kappa|\equiv \kappa_1+\kappa_2+\cdots = k$); 
\begin{equation}
(a)^{\al}_\kappa\equiv\prod_{(i,j)\in\kappa}
\left(a-\frac{i-1}{\alpha}+j-1\right)
\label{poch}
\end{equation}
is the {\em generalized Pochhammer symbol}, and $C_\kappa^{\al}(X)$ is the 
Jack function.

The {\em Jack function} 
$C_\kappa^\al(X)=C_\kappa^\al(x_1,x_2,\ldots,x_n)$ is a symmetric, 
homogeneous polynomial of degree $|\kappa|$ in the eigenvalues 
$x_1,x_2,\ldots,x_n$ of $X$ \cite[Rem.\ 2, p.\ 
228]{muirhead}, \cite{stanley89}. For example, when $\alpha=1$, 
$C_\kappa^\al(X)$ becomes the (normalized) {\em Schur function}, and for 
$\alpha=2$, the {\em zonal polynomial}.\footnote{Some authors define the 
hypergeometric function of a matrix argument through the series 
(\ref{hg1}) for $\alpha=1$~\cite[(4.1)]{grossrichards89} or 
$\alpha=2$~\cite[p.\ 258]{muirhead} only. There is no reason for us to 
treat the different $\alpha$'s separately (see also \cite{dumitriuthesis, 
dumitriuedelman02, forrester} for the uniform treatment of the different 
$\alpha$'s in other settings).} There are several normalizations of the 
Jack function which are scalar multiples of one another: $C_\kappa^\al(X)$ 
is normalized so that $\sum_{\kappa\vdash 
k} C_\kappa^{\al}(X)=(\tr\,X)^k$; in Section \ref{sec_hg} 
we express (\ref{hg1}) in terms of the Jack 
function $J_\kappa^\al(X)$, which is normalized  
so that the coefficient of $x_1x_2\cdots x_{|\kappa|}$  
is $(|\kappa|)!$. The functions $C_\kappa^\al(X)$ and $J_\kappa^\al(X)$ 
can be defined recursively, e.g.,
\begin{equation}
J_\kappa^{\al}(x_1,x_2,\ldots,x_n)=\sum_{\mu\le\kappa}J_\mu^{\al}(x_1,x_2,\ldots,x_{n-1})\cdot 
x_n^{|\kappa/\mu|}\cdot \beta_{\kappa\mu},
\label{recur}
\end{equation}
where $\beta_{\kappa\mu}$ is a rational function of $\alpha$ (see Section 
\ref{sec_hg} for details). The relationship (\ref{recur}) becomes 
key in achieving efficiency in our algorithms.

The Jack functions $C_\kappa^\al(X)$ and $J_\kappa^\al(X)$, and in turn 
the 
hypergeometric function of a matrix argument, depend only on the 
eigenvalues $x_1,x_2,\ldots,x_n$ of $X$. Many authors, however, have found 
the matrix notation in (\ref{hg1}), and the use of a matrix argument to be 
more convenient. We follow the same practice. 

The hypergeometric function of a matrix argument is {\em scalar-valued}, 
which is a major distinction from other functions of a matrix argument 
(e.g., the matrix exponential), which are {\em matrix-valued}. The 
hypergeometric function of a matrix argument generalizes the classical 
hypergeometric function to which it reduces for $n=1$. In general, 
however, there is no explicit relationship between these two functions for 
$n\ge 2$.

We approximate the series (\ref{hg1}) by computing its truncation for 
$|\kappa|\le m$:
\begin{equation}
_p^m F_q^{\al}(a_1,\ldots,a_p;b_1,\ldots,b_q;X)
\equiv \sum_{k=0}^m \sum_{\kappa\vdash k}
\frac{  (a_1)_\kappa^{\al}\cdots(a_p)_\kappa^{\al}}
     {k!(b_1)_\kappa^{\al}\cdots(b_q)_\kappa^{\al}} \cdot
C_\kappa^{\al}(X).
\label{hg2}
\end{equation}

The series~(\ref{hg1}) converges for any $X$ when $p\ge q$; it converges 
if $\max_i|x_i|<1$ when $p=q+1$, and diverges when $p>q+1$, unless it 
terminates~\cite[p.~258]{muirhead}. When it converges, its $\kappa$-term 
converges to zero as $|\kappa|\rightarrow\infty$. In these cases 
(\ref{hg2}) is a good approximation to (\ref{hg1}) for a large enough 
$m$.

The computational difficulties in evaluating (\ref{hg2}) are:
\begin{itemize}

\item[(A)] the series (\ref{hg1}) converges slowly in many 
cases~\cite{muirhead78}; thus a rather large $m$ may be needed before  
(\ref{hg2}) becomes a good approximation to (\ref{hg1});

\item[(B)] the number of terms in (\ref{hg2}) (i.e., the number of 
partitions $|\kappa|\le m$) grows, roughly, as $O(e^{\sqrt{m}})$ (see 
Section \ref{sec_complexity});

\item[(C)] the straightforward evaluation of a single Jack function, 
$C_\kappa^\al(X)$ for $|\kappa|=m$, has 
complexity that grows as $O(n^m)$~\cite{demmelkoevschur}.

\end{itemize}

While there is little we can do about (A) (which is also a major problem 
even in the univariate ($n=1$) case \cite{muller01}), or (B), 
our major 
contribution is in improving (C), the cost of evaluating the Jack 
function. We exploit the combinatorial properties of the Pochhammer symbol 
and the Jack function to only {\em update} the $\kappa$-term in 
(\ref{hg2}) from the $\mu$-terms, $\mu\le\kappa$. As a result the 
complexity of our main algorithm for computing~(\ref{hg2}), 
Algorithm~\ref{alg_hg}, is only linear in the size $n$ of the matrix 
argument $X$, exponentially faster than the previous best 
algorithm~\cite{grs00} (see 
Sections \ref{sec_background}, \ref{sec_complexity} 
and \ref{sec_numexp_complexity} for details). In the special case when 
$X$ is a multiple of the identity, we present an even faster algorithm, 
Algorithm~\ref{alg_hgi}, whose complexity is independent of $n$.

A number of interesting problems remain open. Among these are:
\begin{itemize}
\item detecting convergence;
\item selecting the optimal value of $m$ in (\ref{hg2}) for a 
desired accuracy;
\item selecting the optimal truncation of the series~(\ref{hg1}).
\end{itemize}
We do not believe that a uniform answer to these problems exists for 
every $\alpha$ and every $p$ and $q$. Therefore, we leave the choice of 
$m$ and an appropriate truncation to the user. We elaborate more on these 
open problems in Section \ref{sec_concl}.

With minimal changes our algorithms can approximate the
hypergeometric
function of {\em two} matrix arguments
\begin{equation}
\pFq(a_{1:p};b_{1:q};X;Y)
\equiv\sum_{k=0}^\infty \sum_{\kappa\vdash k}
\frac{(a_1)_\kappa^{\al}\cdots(a_p)_\kappa^{\al}}
   {k!(b_1)_\kappa^{\al}\cdots(b_q)_\kappa^{\al}} \cdot
\frac
{C_\kappa^{\al}(X)\cdot
C_\kappa^{\al}(Y)}
{C_\kappa^{\al}(I)},
\label{hgXY}
\end{equation}
and more generally functions of the form
\begin{equation}
G(X)=\sum_{k=0}^\infty \sum_{\kappa\vdash k} a_\kappa
C_\kappa^{\al}(X),
\label{GX}
\end{equation}
for arbitrary coefficients $a_\kappa$ at a
similar computational cost (see, e.g., (\ref{trA}) 
in subsection \ref{sec_numexp_2}).

In (\ref{hgXY}) and throughout this paper, we denote a vector 
$(z_1,\ldots,z_t)$ as $z_{1:t}$.

This paper is organized as follows.  We survey previous algorithms for 
computing the hypergeometric function of a matrix argument in Section 
\ref{sec_background}.  In Section \ref{sec_hg} we describe our approach in 
computing the truncation (\ref{hg2}).  We present our new algorithms in 
Section \ref{sec_algs}, and analyze their complexity in Section 
\ref{sec_complexity}. We present numerical experiments in Section 
\ref{sec_numexp}. Finally, we draw conclusions and present open problems 
in Section \ref{sec_concl}.

\section{Previous Algorithms}

\label{sec_background}
Butler and Wood~\cite{butlerwood02} used Laplace approximations to compute 
the integral representations~\cite[Thm.\ 7.4.2, p.\ 264]{muirhead}:
\begin{eqnarray*}
_1F_1^{\mbox{\tiny$(2)$}}(a;c;X)
&=&
\frac{\Gamma_n^{\mbox{\tiny$(2)$}}(c)}{\Gamma_n^{\mbox{\tiny$(2)$}}(a)\Gamma_n^{\mbox{\tiny$(2)$}}(c-a)}
\int_{0<Y<I} \!\!\!\!\!e^{\mbox{tr}(XY)}(\det Y)^{a-\frac{n+1}{2}}
\\
&&
\phantom
{
\frac{\Gamma_n^{\mbox{\tiny$(2)$}}(c)}{\Gamma_n^{\mbox{\tiny$(2)$}}(a)\Gamma_n^{\mbox{\tiny$(2)$}}(c-a)}
}
\times \det(I-Y)^{c-a-\frac{n+1}{2}}(dY),
\end{eqnarray*}
valid for real symmetric $X$, 
$\Re(a)>\frac{n-1}{2}$,
$\Re(c)>\frac{n-1}{2}$, and
$\Re(c-a)>\frac{n-1}{2}$; and
\begin{eqnarray*}
_2F_1^{\mbox{\tiny$(2)$}}(a,b;c;X)
&=&
\frac{\Gamma_n^{\mbox{\tiny$(2)$}}(c)}{\Gamma_n^{\mbox{\tiny$(2)$}}(a)\Gamma_n^{\mbox{\tiny$(2)$}}(c-a)}
\int_{0<Y<I} \det(I-XY)^{-b}
\\
&&\hskip0.3in
\times (\det Y)^{a-\frac{n+1}{2}} \det(I-Y)^{c-a-\frac{n+1}{2}}(dY),
\end{eqnarray*}
valid for 
 $\Re(X)<I$,
$\Re(a)>\frac{n-1}{2}$, and
$\Re(c-a)>\frac{n-1}{2}.$

This approach, however, is restricted to the cases $p=1$ or $2$, $q=1$, 
and $\alpha=2$.

Guti\'errez, Rodriguez, and S\'aez presented in~\cite{grs00} (see
also~\cite{grs00soft} for the implementation)
an algorithm for computing the truncation~(\ref{hg2}) for
$\alpha=2$ (then the Jack functions are called {\em zonal polynomials}).
For every
$k=1,2,\ldots, m$, the authors form the upper triangular {\em transition
matrix}~\cite[p.~99]{macdonald} $K$ (indexed by all partitions of
$k$) between the monomial symmetric functions
$(m_\kappa)_{\kappa\vdash k}$ and the zonal polynomials
$(C_\kappa)_{\kappa\vdash k}$.
Then for every partition $\kappa=(\kappa_1,\kappa_2,\ldots)\vdash k$ they
compute
$$
m_\kappa=\sum_\mu x_1^{\mu_1} x_2^{\mu_2}\cdots
$$
(where $\mu$ ranges over all {\em distinct} permutations of
$(\kappa_1,\kappa_2,\ldots)$~\cite[p.~289]{EC2}), and form the product
$(C_\kappa)_{\kappa\vdash k} = K\cdot (m_\kappa)_{\kappa\vdash k}$.
Computing the vector $(m_\kappa)_{\kappa\vdash m}$ alone costs
$m\binom{n+m-1}{m}=O(n^m)$ since every term in every $m_\kappa$ is of
degree $m$, and for every nonstrictly increasing sequence of $m$ numbers
from the set $\{1,2,\ldots,n\}$ we obtain a distinct term in some
$m_\kappa$.  The overall cost is thus at least exponential
($O(n^m)$),
which explains the authors' observation:
\begin{quotation}
We spent about 8 days to obtain the 627 zonal polynomials of degree
20 with a 350 MHz Pentium II processor.
\end{quotation}

In contrast, our Algorithm~\ref{alg_hg} takes less than a hundredth of a
second to do the same. Its complexity is only linear in $n$ and
subexponential in $m$ (see Section \ref{sec_complexity}).

\section{Our Approach}
\label{sec_hg}

We make the evaluation of $\mpFq(a_{1:p};b_{1:q};X)$ efficient by 
only {\em updating} the $\kappa$-term from the $\mu$-terms, 
$\mu\le\kappa$, 
instead of computing it from scratch.

We first express $\mpFq(a_{1:p};b_{1:q};X)$ in
terms of the Jack function $J_\kappa^{\al}(X)$, which is normalized
so that the coefficient of $x_1x_2\ldots x_{|\kappa|}$ in
$J_\kappa^{\al}(X)=J_\kappa^{\al}(x_1,x_2,\ldots,x_n)$ equals
$(|\kappa|)!$~\cite[Thm.\ 1.1]{stanley89}. 
The Jack functions $C_\kappa^{\al}(X)$ and
$J_\kappa^{\al}(X)$ are related as:
\begin{equation}
C_\kappa^{\al}(X)
=\frac{\alpha^{|\kappa|} \cdot (|\kappa|)!}{j_\kappa}J_\kappa^{\al}(X),
\label{ckappajkappa}
\end{equation}
where 
\begin{equation}
j_\kappa=\prod_{(i,j)\in\kappa}h_{*}^\kappa(i,j)h_\kappa^{*}(i,j),
\label{jkappa}
\end{equation}
and $
h_\kappa^{*}(i,j)\equiv\kappa_j'-i+\alpha(\kappa_i-j+1)
$ and $h_{*}^\kappa(i,j)\equiv\kappa_j'-i+1+\alpha(\kappa_i-j)$
are the {\em upper and
lower hook lengths} at $(i,j)\in\kappa$, respectively.

Denote
\begin{equation}
Q_\kappa\equiv\frac{\alpha^{|\kappa|}
(a_1)_\kappa^{\al}\ldots(a_p)_\kappa^{\al}}
{j_\kappa(b_1)_\kappa^{\al}\ldots(b_q)_\kappa^{\al}}.
\label{Qkappa}
\end{equation}
Since $J_\kappa^{\al}(X)=0$ when $\kappa_{n+1}>0$, we need to sum 
only over partitions $\kappa$ with at most $n$ parts:
\begin{eqnarray}
^m_p F^{\al}_q(a_{1:p};
b_{1:q};X)
&=&\sum_{|\kappa|\le m,\, \kappa_{n+1}=0} Q_\kappa
J_\kappa^{\al}(X).
\label{hg3}
\end{eqnarray}

When computing~(\ref{hg3}), we recursively generate all partitions
$|\kappa|\le m$ in such a way that consecutively generated partitions
differ in only one part. Therefore it is convenient to introduce the 
notation
$$
\kappa_{(i)}\equiv(\kappa_1,\ldots,\kappa_{i-1},\kappa_i-1,
\kappa_{i+1},\ldots)
$$
for any partition $\kappa$ such that $\kappa_i>\kappa_{i+1}$.

In the following subsections we derive formulas for updating
the $\kappa$-term in~(\ref{hg3}) from the $\mu$-terms, $\mu\le\kappa$.

\subsection{Updating the Coefficients $Q_\kappa$}

We update $Q_\kappa$ from $Q_{\kappa_{(i)}}$ using the following
lemma.
\begin{lemma}
\begin{equation}
\frac{Q_\kappa}{Q_{\kappa_{(i)}}}=
\frac
{\prod_{j=1}^p(a_j+c)}
{\prod_{j=1}^q (b_j+c)}
\cdot
\prod_{j=1}^{\kappa_i-1}
\frac
{(g_j-\alpha)e_j}
{g_j(e_j+\alpha)}
\cdot
\prod_{j=1}^{i-1}
\frac
{l_j-f_j}
{l_j+h_j},
\label{Qkappaupdate}
\end{equation}
where $c=-\frac{i-1}{\alpha}+\kappa_i-1,\, d=\kappa_i\alpha-i,\,
e_j=d-j\alpha+\kappa_j',\;g_j=e_j+1,\;f_j=\kappa_j\alpha-j-d,\,
h_j=f_j+\alpha$, and
$l_j=h_jf_j$.
\label{lemma_Qkappaupdate}
\end{lemma}

\begin{proof}
From~(\ref{poch}),
$(a)_\kappa^{\al}=(a)_{\kappa_{(i)}}^{\al}\cdot
(a-(i-1)/\alpha+\kappa_i-1)$, and from~(\ref{jkappa}),
\begin{equation}
\frac
{j_{\kappa_{(i)}}}
{j_\kappa}
= \frac{1}{\alpha} \cdot
\prod_{j=1}^{\kappa_i-1}
\frac
{h_{*}^{\kappa_{(i)}}(i,j)\cdot h_{\kappa_{(i)}}^{*}(i,j)}
{h_{*}^\kappa(i,j)\cdot h_\kappa^{*}(i,j)}
\cdot
\prod_{j=1}^{i-1}
\frac
{h_{*}^{\kappa_{(i)}}(j,\kappa_i)\cdot
h_{\kappa_{(i)}}^{*}(j,\kappa_i)}
{h_{*}^\kappa(j,\kappa_i)\cdot h_\kappa^{*}(j,\kappa_i)},
\label{Qkappaupdate2}
\end{equation}
which along with~(\ref{Qkappa}) imply~(\ref{Qkappaupdate}).
\end{proof}

The seemingly complicated notation of Lemma \ref{lemma_Qkappaupdate} is
needed in order to minimize the number of arithmetic operations needed to
update $Q_\kappa$. A straightforward evaluation
of~(\ref{Qkappa}) costs $6|\kappa|(2+p+q)$; 
in contrast,~(\ref{Qkappaupdate}) costs only $2(p+q)+10\kappa_i+9i-11$
arithmetic operations.

\subsection{Updating the Jack Function}
\label{sec_updJack}

When $\kappa=(0)$, $J_{(0)}^\al(x_1,\ldots,x_n)=1$. For 
$\kappa>(0)$, we update 
$J_\kappa^\al(x_1,\ldots,x_n)$ from 
$J_\mu^\al(x_1,\ldots,x_r),\,\mu\le\kappa,\,r\le n.$

When $X$ is a multiple of the identity we have an easy special 
case~\cite[Thm.~5.4]{stanley89}: 
$$
J_\kappa^{\al}(x I)=
x^{|\kappa|}
\prod_{(i,j)\in \kappa}
(n-(i-1)+\alpha(j-1)).
$$
Therefore we can update $J_\kappa^{\al}(xI)$ from
$J_{\kappa_{(i)}}^{\al}(xI)$ as
\begin{equation}
J_\kappa^{\al}(x I)=J_{\kappa_{(i)}}^{\al}(xI)\cdot
x\cdot (n-i+1+\alpha(\kappa_i-1)).
\label{Jkappaupdate}
\end{equation}

In the general case, we update the Jack function using the
identity (see, e.g., \cite[Prop.\ 4.2]{stanley89}):
\begin{equation}
J_\kappa^{\al}(x_1,x_2,\ldots,x_n)=\sum_{\mu\le\kappa}
J_\mu^{\al}(x_1,x_2,\ldots,x_{n-1})
x_n^{|\kappa/\mu|}\beta_{\kappa\mu},
\label{jackmain}
\end{equation}
where the summation is over all $\mu\le\kappa$ such that
$\kappa/\mu$ is a horizontal strip, and
\begin{equation}
\beta_{\kappa\mu}\equiv\frac{
 \prod_{(i,j)\in \kappa} B_{\kappa\mu}^\kappa(i,j)
}{
\prod_{(i,j)\in \mu} B_{\kappa\mu}^\mu(i,j)
},
\mbox{\hskip0.2in where\hskip0.2in}
B_{\kappa\mu}^\nu(i,j)\equiv\left\{
\begin{array}{ll}
h_\nu^{*}(i,j), & \mbox{if }\kappa_j'=\mu_j'; \\
h_{*}^\nu(i,j), & \mbox{otherwise.}
\end{array}
\right.
\label{betakappamu}
\end{equation}
The skew partition $\kappa/\mu$ is a {\em horizontal strip} when
$\kappa_1\ge\mu_1\ge\kappa_2\ge\mu_2\ge\ldots$~\cite[p.~339]{EC2}.

We borrow the idea for updating the Jack function from
\cite{demmelkoevschur}, but make two
important improvements. We only {\em update} the coefficients
$\beta_{\kappa\mu}$,
and store the precomputed Jack functions much more efficiently than
in~\cite{demmelkoevschur}.

The coefficients $\beta_{\kappa\mu}$ are readily computable 
using~(\ref{betakappamu}) at the cost of $6(|\kappa|+|\mu|)$ arithmetic 
operations.  The following lemma allows us to start with 
$\beta_{\kappa\kappa}=1$ and {\em update} 
$\beta_{\kappa\mu_{(k)}}$ from $\beta_{\kappa\mu}$ at the cost of only 
$12k+6\mu_k-7$ arithmetic operations.

\begin{lemma}
Let $\kappa, \mu$, and $\nu=\mu_{(k)}$ be partitions such that
$\kappa/\mu$ and $\kappa/\nu$ are horizontal strips,
 and $\kappa_r'=\mu_r'$ for $0\le r\le k-1$. Then
\begin{equation}
\frac{\beta_{\kappa\nu}}{\beta_{\kappa\mu}}=
\alpha\cdot
\prod_{r=1}^k
\frac{u_r}{u_r+\alpha'}\cdot
\prod_{r=1}^{k-1}\frac{v_r+\alpha}{v_r}
\cdot \prod_{r=1}^{\mu_k-1} \frac{w_r+\alpha}{w_r},
\label{betakappamuupdate}
\end{equation}
where $\alpha'=\alpha-1,\, t=k-\alpha\mu_k,\, 
q=t+1, \, u_r=q-r+\alpha\kappa_r,
v_r=t-r+\alpha\mu_r$, and $w_r=\mu_r'-t-\alpha r$.
\end{lemma}

\begin{proof}
Denote $l\equiv\mu_k$. We have $\nu=\mu$, except for $\nu_k=\mu_k-1$.
Since $\kappa/\mu$ and $\kappa/\nu$ are horizontal strips,
$\kappa_l'=\mu_l'$ and $\kappa_l'\ne\nu_l'=\mu_l'-1$. Then
\begin{eqnarray}
\frac{\beta_{\kappa\nu}}{\beta_{\kappa\mu}}
&=&
\frac
{\prod_\mu     B_{\kappa\mu}^\mu}
{\prod_\kappa B_{\kappa\mu}^\kappa}
\cdot
\frac
{\prod_\kappa B_{\kappa\nu}^\kappa}
{\prod_{\nu} B_{\kappa\nu}^{\nu}}
\nonumber \\
&=&
\prod_\kappa
\frac
{B_{\kappa\nu}^\kappa}
{B_{\kappa\mu}^\kappa}
\cdot
\prod_{\nu}
\frac
{B_{\kappa\mu}^\mu}
{B_{\kappa\nu}^{\nu}}
\cdot B_{\kappa\mu}^\mu(k,l)
\nonumber \\
&=&
\prod_\kappa
\frac
{B_{\kappa\nu}^\kappa}
{B_{\kappa\mu}^\kappa}
\cdot
\prod_{\nu}
\frac
{B_{\kappa\mu}^\mu}
{B_{\kappa\nu}^{\nu}}
\cdot \alpha.
\label{betatrans}
\end{eqnarray}
We transform the first term of~(\ref{betatrans})
by observing that
$B_{\kappa\mu}^\kappa(i,j)=B_{\kappa\nu}^\kappa(i,j)$, except
for $j=l$:
$$
\prod_\kappa
\frac
{B_{\kappa\nu}^\kappa}
{B_{\kappa\mu}^\kappa}
=
\prod_{r=1}^k
\frac
{B_{\kappa\nu}^\kappa(r,l)}
{B_{\kappa\mu}^\kappa(r,l)}
=
\prod_{r=1}^k
\frac{h^\kappa_*(r,l)}{h_\kappa^*(r,l)}
=
\prod_{r=1}^k
\frac{
k-r+1+\alpha(\kappa_r-l)
}{
k-r+  \alpha(\kappa_r-l+1)
}.
$$
To simplify the second term of~(\ref{betatrans}), we observe that
$B_{\kappa\mu}^\mu(i,j)=B_{\kappa\nu}^\nu(i,j)$,
except for $i=k$ and $j=l$. Also $\kappa_r'=\mu_r'=\nu_r'$ for $0\le r\le
k-1$, and $\kappa_k'=\mu_k'\ne\mu_k'-1=\nu_k'$. Therefore
\begin{eqnarray*}
\prod_{\nu}
\frac
{B_{\kappa\mu}^\mu}
{B_{\kappa\nu}^{\nu}}
&=&
\prod_{r=1}^k
\frac
{B_{\kappa\mu}^\mu(r,l)}
{B_{\kappa\nu}^{\nu}(r,l)}
\cdot
\prod_{r=1}^l
\frac
{B_{\kappa\mu}^\mu(k,r)}
{B_{\kappa\nu}^{\nu}(k,r)}
\\
&=&
\prod_{r=1}^{l-1}
\frac
{
h_\mu^*(k,r)
}
{
h_{\nu}^*(k,r)
}
\cdot
\prod_{r=1}^{k-1}
\frac
{
h_\mu^*(r,l)
}
{
h^{\nu}_*(r,l)
}
\\
&=&
\prod_{r=1}^{l-1}
\frac
{\mu_r'-k+\alpha(\mu_k-r+1)}
{\mu_r'-k+\alpha(\mu_k-r)}
\cdot
\prod_{r=1}^{k-1}
\frac
{k-r+\alpha(\mu_r-l+1)}
{k-r+\alpha(\mu_r-l)},
\end{eqnarray*}
which yields (\ref{betakappamuupdate}).
 \end{proof}

\section{Algorithms for Efficient Evaluation of
$_p^mF_q^{(\alpha)}(a_{1:p};b_{1:q};X)$}
\label{sec_algs}

Algorithm~\ref{alg_hgi} computes
$\mpFq(a_{1:p};b_{1:q};X)$ in the easy special case when
$X$ is a multiple of the identity. Algorithm~\ref{alg_hg}
handles the general case.

Both algorithms recursively generate all partitions $|\kappa|\le m$ with at
most $n$ parts by allowing $\kappa_1$ to take all values $1,2,\ldots,m$
and, independently, $\kappa_i, i=1,2,\ldots,n$, to take all values
$1,2,\ldots,\kappa_{i-1}$, subject to the restriction $|\kappa|\le m$.
The coefficients $Q_\kappa$ are updated using ~(\ref{Qkappaupdate}). The
Jack functions are updated using~(\ref{Jkappaupdate}) or~(\ref{jackmain})
as appropriate.

\subsection{The Case when $X$ is a Multiple of the Identity}

\begin{algorithm}[Hypergeometric Function, $X=xI$]
The following algorithm computes $\mpFq(a_{1:p};b_{1:q};xI)$.
 The variables $x, n, \alpha, s,a,b$, and $\kappa$ are global.

\smallskip

\noindent\verb+ function +$s=$ \verb+hgi+$(m,\alpha,a,b,n,x)$

\noindent\verb+   +$s=1$

\noindent\verb+   summation+$(1,1,m)$

\smallskip

\noindent\verb+ function summation+$(i,z,j)$

\noindent\verb+   for +$\kappa_i=1:
                     \min(\kappa_{i-1},j)$
                      \hskip1.8in ({\em defaults to $j$ for $i=1$})

\ignore{

\noindent\verb+     +$\kappa'_{\kappa_i}=\kappa'_{\kappa_i}+1$

\noindent\verb+        +$z=z\cdot \verb+prod+(a+c)/\verb+prod+(b+c)$

\noindent\verb+        +$d=\kappa_i\alpha-i$

\noindent\verb+        for +$j=1:\kappa_i-1$

\noindent\verb+          +$e=d-j\alpha+\kappa'_j;\;\;g=e+1$

\noindent\verb+          +$z=z(g-\alpha)e/(g(e+\alpha))$

\noindent\verb+        for +$j=1:i-1$

\noindent\verb+          +$f=\kappa_j\alpha-j-d;\;\;h=f+\alpha;\;\;l=hf$

\noindent\verb+          +$z=z(l-f)/(l+h)$
}


\noindent\verb+     +$z=z\cdot x\cdot (n-i+1+\alpha(\kappa_i-1))\cdot T$
                        \hskip0.1in
                    ({\em where $T$ is the right hand side
                     of} (\ref{Qkappaupdate}))

\noindent\verb+     +$s=s+z$

\noindent\verb+     if +$(j>\kappa_i)$\verb+ and +$(i<n)$\verb+ then+

\noindent\verb+       summation+$(i+1,z,j-\kappa_i)$

\noindent\verb+     endif+

\noindent\verb+   endfor+

\noindent\verb+   +$\kappa_i=0$

\medskip

\label{alg_hgi}
\end{algorithm}

In Algorithm~\ref{alg_hgi} the variable $z$ equals the $\kappa$-term
in~(\ref{hg3}); it is updated using~(\ref{Qkappaupdate})
and~(\ref{Jkappaupdate}). The parameter $j$ in \verb+summation+ equals
$m-|\kappa|$.

The hypergeometric function of two matrix arguments~(\ref{hgXY})
is in this case $$\pFq(a_{1:p};b_{1:q}; xI;
yI)=~\!\pFq(a_{1:p};b_{1:q};  xyI).$$

\subsection{The General Case} We use the identity~(\ref{jackmain}) to 
update the Jack function, but in order to use it efficiently, we need to 
store and reuse the Jack functions computed earlier.

Therefore we index all partitions $|\kappa|\le m$ with at most $n$ parts 
($\kappa_{n+1}=0$) by linearizing the
$m$-tree that they form (each node
$\kappa=(\kappa_1,\kappa_2,\ldots,\kappa_k)$ has at most $m$ children
$(\kappa_1,\kappa_2,\ldots,\kappa_k,\kappa_{k+1}),$
$\kappa_{k+1}=1,2,\ldots,\min(\kappa_k,m-|\kappa|),\;k<n$).
In other words, if $P_{mn}=\#\{\kappa\vert \; |\kappa|\le 
m,\;\kappa_{n+1}=0\}$, we assign a 
distinct integer index
$N_\kappa\in\{0,1,\ldots,P_{mn}\}$ to every such partition $\kappa$.
We start by assigning the indexes $0,1,\ldots,m$ to
partitions with one part: $N_{(i)}\equiv i$,  $i=0,1,2,\ldots,m$. Then,
recursively, once an index $N_\kappa$ has been assigned to a partition
$\kappa=(\kappa_1,\ldots,\kappa_k)$ with $k$ parts, we assign $m-|\kappa|$
consecutive unassigned indexes to the partitions
$(\kappa_1,\ldots,\kappa_k,\kappa_{k+1})$,
$\kappa_{k+1}=1,2,\ldots,m-|\kappa|$.
We record the tree structure in an array $D$ such that
$D(N_{(\kappa_1,\ldots,\kappa_k)})= N_{(\kappa_1,\ldots,\kappa_k,1)}$. Now
given $\kappa$, we can compute $N_\kappa$ by starting with
$N_{(\kappa_1)}=\kappa_1$ and using the recurrence
\begin{equation}
N_{(\kappa_1,\ldots,\kappa_i)}=
D(N_{(\kappa_1,\ldots,\kappa_{i-1})})+\kappa_i-1.  \label{remU}
\end{equation}

We store every computed $J_\kappa^{\al}(x_{1:i})$ ($|\kappa|\le m$, 
$\kappa_{n+1}=0$, $i=1,2,\ldots, n$) in the $(N_\kappa,i)$th entry of an 
$P_{mn}\times n$ array, which we call ``$J$'' in Algorithm~\ref{alg_hg} 
below.

We compute the value of $P_{mn}$ as follows. Let $p_k(i)$ be the number
of partitions of $i$ with exactly $k$ parts.
The $p_k(i), i=1,2,\ldots,m, k=1,2,\ldots,\min(m,n)$, are computed using 
the recurrence $p_k(i)=p_{k-1}(i-1)+p_k(i-k)$~\cite[p.~28]{EC1}. Then
\begin{equation}
P_{mn}=\sum_{i=1}^m\sum_{k=1}^{\min(m,n)} p_k(i).
\label{remarkPm}
\end{equation}

Next, we present our main algorithm.

\begin{algorithm}[Hypergeometric Function]
The following algorithm computes
$\mpFq(a_{1:p};b_{1:q};X)$, where $X=\mbox{diag}(x)$.
The variables $m,n,a,b,\alpha,x,J,D,H,\kappa$, $N_\kappa,\mu,N_\mu$, and
$s$ are global.

\smallskip

\noindent\verb+ function +$s=\verb+hg+(m,\alpha,a,b,x)$

\noindent\verb+   +Compute $P_{mn}$ using~(\ref{remarkPm})

\noindent\verb+   +$n=\verb+length+(x);\;\;H=m+1;\;\;s=1;\;\;
                    D=\verb+zeros+(P_{mn},1)$

\noindent\verb+   +$J=\verb+zeros+(P_{mn},n);\;J(1,:)=1$
                                   \hskip0.9in({\em
                   $J_\kappa^{\al}(x_{1:i})$ is stored in
                   $J(N_\kappa,i)$})

\noindent\verb+   summation+$(1,1,1,m)$

\smallskip

\noindent\verb+ function summation+$(i,z,j)$

\noindent\verb+   +$r=\min(\kappa_{i-1},j)$
                     \hskip2.3in ({\em defaults to $j$ for $i=1$})

\noindent\verb+   for +$\kappa_i=1:r$

\noindent\verb+     if +$(\kappa_i=1)\verb+ and +(i>1)$\verb+ then +
                    $D(N_\kappa)=H;\;\; N_\kappa=H;\;\;H=H+r$

\noindent\verb+     else +$N_\kappa=N_\kappa+1$

\noindent\verb+     endif+

\noindent\verb+     +$z=z\cdot T$ \hskip1.65in ({\em where $T$
                     is the right hand side of} (\ref{Qkappaupdate}))

\noindent\verb+     if +$\kappa_1'=1$\verb+ then+

\noindent\verb+       +$J(N_\kappa,1)=x_1 (1+\alpha(\kappa_1-1))
                          \cdot J(N_\kappa-1,1)$

\noindent\verb+     endif+

\noindent\verb+     for +$t=2:n$\verb+ do jack+$(0,1,0,t)$\hskip1.4in
                     ({\em computes $J_\kappa^{\al}(x_{1:t})$})

\noindent\verb+     endfor+

\noindent\verb+     +$s=s+z\cdot J(N_\kappa,n,1)$

\noindent\verb+     if +$(j>\kappa_i)\verb+ and +(i<n)$\verb+ then+

\noindent\verb+       summation+$(i+1,z,j-\kappa_i)$

\noindent\verb+     endif+

\noindent\verb+   endfor+

\smallskip

\noindent\verb+ function jack+$(k,\beta_{\kappa\mu},c,t)$

\noindent\verb+   for +$i=k:\mu'_1$

\noindent\verb+     if +$(k>0)\verb+ and +(\mu_i>\mu_{i+1})$

\noindent\verb+       +$d=N_\mu$

\noindent\verb+       +$\mu_i=\mu_i-1;\;\;$
                    Compute $N_\mu$ using~(\ref{remU})

\noindent\verb+       +Update $\beta_{\kappa\mu}$
using~(\ref{betakappamuupdate})

\noindent\verb+       if +$\mu_i>0$\verb+ then jack+$
                           (i,\beta_{\kappa\mu},c+1,t)$

\noindent\verb+       else +$
                            J(N_\kappa,t)=
                            J(N_\kappa,t)+\beta_{\kappa\mu}
                            \cdot
                            J(N_\mu,t-1)\cdot x_t^{c+1}$

\noindent\verb+       endif+

\noindent\verb+       +$\mu_i=\mu_i+1;\; N_\mu=d$

\noindent\verb+     endif+

\noindent\verb+   endfor+

\noindent\verb+   if +$k=0$\verb+ then +$J(N_\kappa,t)=
                           J(N_\kappa,t)+J(N_\kappa,t-1)$

\noindent\verb+   else +$J(N_\kappa,t)=
                           J(N_\kappa,t)+\beta_{\kappa\mu}
                           \cdot x_t^c$

\noindent\verb+   endif+

\label{alg_hg}
\end{algorithm}

Algorithm~\ref{alg_hg} generates all partitions $|\kappa|\le m$ with at
most $n$ parts analogously to Algorithm~\ref{alg_hgi}. The parameter $z$
in \verb+summation+ is now just $Q_\kappa$.

For every $\kappa$ the function \verb+jack+ recursively generates all
partitions $\mu\le\kappa$ such that $\kappa/\mu$ is a horizontal strip.
The Jack function is computed using~(\ref{jackmain}); $\beta_{\kappa\mu}$
is updated using~(\ref{betakappamuupdate}).  The parameter $c$ equals
$|\kappa|-|\mu|$ at all times.

\subsection{Implementation notes}

In our implementation of
Algorithms~\ref{alg_hgi} and~\ref{alg_hg} in~\cite{koevwebpage} we made a
few fairly straightforward,
but important improvements:
\begin{enumerate}
\item we precompute the values of $x_i^j$ for $i=1,2,\ldots,n$, and
$j=1,2,\ldots,m$;
\item we keep and update the conjugate partitions $\kappa'$ and $\mu'$
along with $\kappa$ and $\mu$, so we never need to recover $\kappa'$ and
$\mu'$ when  computing~(\ref{Qkappaupdate})
and~(\ref{betakappamuupdate});
\item when two sets of arguments ($x_{1:n}$ and $y_{1:n}$) are passed,
the hypergeometric function of two
matrix arguments, (\ref{hgXY}), is computed;
\item if a vector $x=(t_1,\ldots,t_r)$ is passed as a parameter in
Algorithm~\ref{alg_hgi}, the output is also a vector with the values
of $_p^m F^{\al}_q(a_{1:p};b_{1:q}; t_jI)$ for $j=1,2,\ldots,r.$

\end{enumerate}

\section{Complexity Analysis}
\label{sec_complexity}

Algorithm~\ref{alg_hgi} (the case $X=xI$) costs $O(P_{mn})$
arithmetic operations (where again, $P_{mn}$ is the number of partitions 
$\kappa$,
$|\kappa|\le m$ with at most $n$ parts).

\ignore{
It costs at
most $2(p+q)+9m$ arithmetic operations to evaluate~(\ref{Qkappaupdate}),
therefore
the cost of
}

To bound the complexity of Algorithm~\ref{alg_hg} (general case) we
observe that the formula~(\ref{jackmain}) represents the summation of at
most $P_{mn}$ terms (in fact a lot less, but we have been unable
to obtain a better bound than $P_{mn}$ that is easy to work with). We
use~(\ref{jackmain}) to compute the Jack function
$J_\kappa^{\mbox{\tiny($\alpha$)}}(x_{1:t})$ for all $|\kappa|\le m$
and $t\le n$, i.e., $P_{mn}\cdot n$ Jack functions, each of which costs at
most $O(P_{mn})$ arithmetic operations. The overall cost of
Algorithm~\ref{alg_hg} is thus bounded by
$$
O(P_{mn}^2\cdot n).
$$

There is no explicit formula for $P_{mn}$, so we use Ramanujan's
asymptotic formula~\cite[p.~116]{hardy99}
for number of partitions of $m$
$$
p(m)\sim\frac{1}{4m\sqrt{3}} \exp\big(\pi \sqrt{2m/3}\big)
$$
to obtain
$$
P_{mn}\le \sum_{i=1}^mp(i) \sim
O\left(\exp\big(\pi\sqrt{2m/3}\big)\right),\;\;\;\mbox{~~~and~~~} \;\;\;
P_{mn}^2\sim O\left(\exp\big(2\pi\sqrt{2m/3}\big)\right).
$$
Therefore the complexity of Algorithm~\ref{alg_hg} is {\em linear} in
$n$ and subexponential in~$m$.

\section{Numerical Experiments}
\label{sec_numexp}

We performed extensive numerical tests to verify the correctness and 
complexity of our Algorithms~\ref{alg_hgi} and~\ref{alg_hg}.  We compared 
the output of our algorithms for $_0F^{\al}_0$ and $_1F^{\al}_1$ with 
explicit expressions (subsection \ref{sec_explicit}).  We also compared 
the probability distributions of the eigenvalues of certain random 
matrices (which are expressed as hypergeometric functions) against the 
results of Monte--Carlo experiments (subsection \ref{sec_numexp_2}). 
Finally, we present performance results in 
subsection~\ref{sec_numexp_complexity}.

\subsection{Explicit Expressions}
\label{sec_explicit}

We compared the output of Algorithms~\ref{alg_hgi} and~\ref{alg_hg}
against the expressions~\cite[p.\ 262]{muirhead}:
\begin{eqnarray}
\phantom{^(}_0F_0^{\al}(X)&=&e^{\tr\, X};
\label{0F0_s}
\\
\phantom{^(}_1F_0^{\al}(a;X)&=&\det (I-X)^{-a},
\label{1F0_s}
\end{eqnarray}
for $n=10$ and 
random uniformly distributed values $x_i\in[0,\frac{1}{2}]$, 
$i=1,2,\ldots,n$. For $m=30$ it took less than one second per test. The 
results agreed to at least $12$ decimal digits with~(\ref{0F0_s}), and at 
least $8$ decimal digits with~(\ref{1F0_s}), reflecting the slower 
convergence of~(\ref{1F0_s}).

\subsection{Eigenvalue Statistics}
\label{sec_numexp_2}
We tested Algorithms~\ref{alg_hgi} and~\ref{alg_hg} against
the eigenvalue statistics of the $\beta$-Laguerre and
Wishart matrices.

The $n\times n$ {\em $\beta$-Laguerre matrix}
of parameter $a$ is defined as $L=BB^T$, where
$$
B=\bmat{cccc}
\chi_{2a} \\
\chi_{\beta(n-1)}  & \chi_{2a-\beta}\\
& \ddots & \ddots\\
&& \chi_\beta
&\chi_{2a-\beta(n-1)}
\emat, \;\;\;
a>\frac{\beta}{2}(n-1).
$$

The $n\times n$ {\em Wishart matrix} with $l$ degrees
of freedom ($l>n$) and covariance matrix $\Sigma$, $A=W_n(l,\Sigma)$,
is defined as $A=\Sigma^{1/2}\cdot Z^T \cdot Z\cdot \Sigma^{1/2}$, where
$Z_{ij}=N(0,1)$ for $i=1,2,\ldots,l$; $j=1,2,\ldots,n$.

The eigenvalue distributions of $A$ and $L$ are the same when $\beta=1$
or $2$, $a=\frac{l}{2}$, and $\Sigma=I$.

\begin{figure}[t]
\centerline{\resizebox{5in}{3.95in}{
\includegraphics*{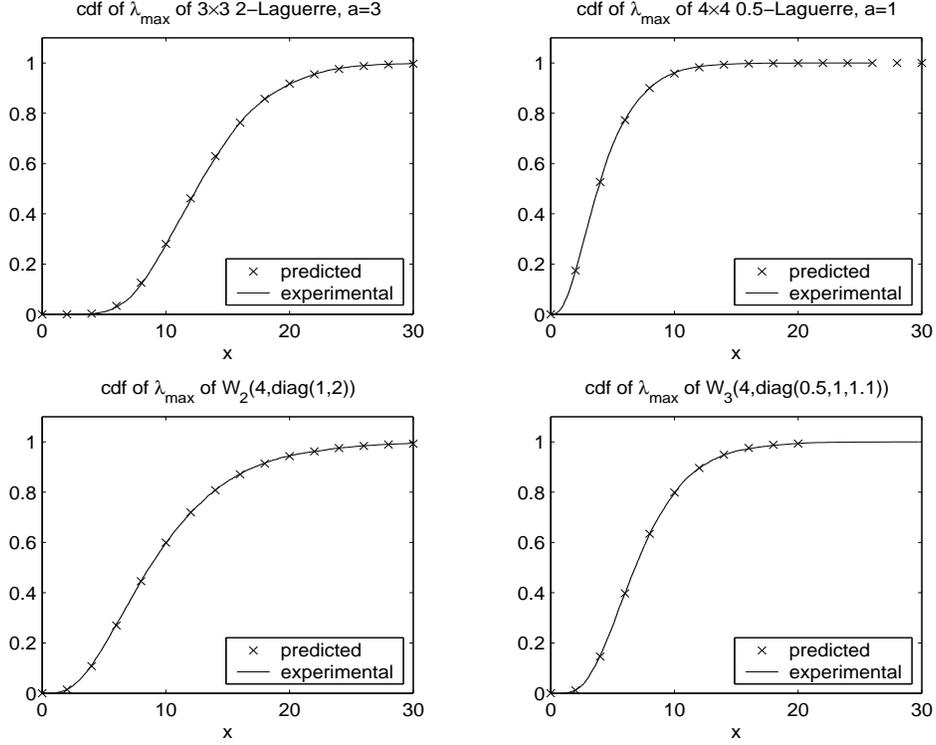}
}}
\caption{
The c.d.f.~of $\lambda_{\max}$ of the $\beta$-Laguerre matrix
(top plots) and the Wishart matrix (bottom plots).}
\label{fig_example1}
\end{figure}

The cumulative distribution functions of the largest eigenvalues,
$\lambda_L$ and $\lambda_A$, of $L$ and
$A$, are
\begin{eqnarray}
P(\lambda_L<x)
&=&\frac{\Gamma_n^{\al}\big(\frac{n-1}{\alpha}+1\big)}
     {\Gamma_n^{\al}\big(a+\frac{n-1}{\alpha}+1\big)}
     \left(\frac{x}{2}\right)^{an}\!\!
     \phantom{^2}_1F_1^{\al}\big(a;a+{\textstyle
       \frac{n-1}{\alpha}+1;-\frac{1}{2}x I}
     \big),
\nonumber 
\\
P(\lambda_A<x)
&=&\frac{\Gamma_n^{\mbox{\tiny$(2)$}}\big({\textstyle\frac{n+1}{2}}\big)}
     {\Gamma_n^{\mbox{\tiny$(2)$}}\big({\textstyle\frac{l+n+1}{2}}\big)}
     \det\big({\textstyle\frac{1}{2}}x\Sigma^{-1}\big)^{l/2}
     \!\!
     \phantom{^2}_1F_1^{\mbox{\tiny$(2)$}}\big(
     {\textstyle\frac{l}{2}};{\textstyle\frac{n+l+1}{2}};
     -{\textstyle\frac{1}{2}} x\Sigma^{-1}\big),
\nonumber 
\end{eqnarray}
respectively \cite[Thm.~10.2.1, p.~147]{dumitriuthesis}, \cite[Thm.~9.7.1, 
p.~420]{muirhead}, where $\alpha=2/\beta$, and $\Gamma_n^{\al}$ is the {\em 
multivariate Gamma function of parameter $\alpha$}:
\begin{equation}
\Gamma_n^{\al}(c)\equiv\pi^{\frac{n(n-1)}{2\alpha}} \prod_{i=1}^n
\Gamma\left(c-\frac{i-1}{\alpha}\right)
\mbox{~~~~for~~}\Re(c)>\frac{n-1}{\alpha}.
\label{gamman}
\end{equation}
We use the {\em Kummer relation}~\cite[Thm.\ 7.4.3, p.\ 265]{muirhead}
$$
\phantom{^\al}_1F_1^{\al}(a;c;X)
= e^{\tr\, X} \cdot\!\!\!
\phantom{^\al}_1F_1^{\al}(c-a;c;-X)
$$
to obtain the equivalent, but 
numerically more stable expressions
\begin{eqnarray}
P(\lambda_L<x)
\!\!&\!\!=\!\!&\!\!
\frac{\Gamma_n^{\al}\big(\frac{n-1}{\alpha}+1\big)}
     {\Gamma_n^{\al}\big(a+\frac{n-1}{\alpha}+1\big)}
     \left(\frac{x}{2}\right)^{an} e^{-\frac{nx}{2}}\!\!
     \phantom{^2}_1F_1^{\al}\big({\textstyle
            \frac{n-1}{\alpha}+1;
            a+\frac{n-1}{\alpha}+1;\frac{1}{2}x I}
     \big);
\nonumber 
\\
P(\lambda_A<x)
\!\!&\!\!=\!\!&\!\!
\frac{\Gamma_n^{\mbox{\tiny$(2)$}}\big({\textstyle\frac{n+1}{2}}\big)}
     {\Gamma_n^{\mbox{\tiny$(2)$}}\big({\textstyle\frac{l+n+1}{2}}\big)}
     \det\big({\textstyle\frac{1}{2}}x\Sigma^{-1}\big)^{l/2}
 e^{\tr(-\frac{x}{2}\Sigma^{-1})}
     \!\!
     \phantom{^2}_1F_1^{\mbox{\tiny$(2)$}}\big(
     {\textstyle\frac{n+1}{2}};{\textstyle\frac{n+l+1}{2}};
     {\textstyle\frac{1}{2}} x\Sigma^{-1}\big),
\nonumber 
\end{eqnarray}
which we plot on Figure~\ref{fig_example1} along with the Monte--Carlo
results from a sample of 10000 random matrices.

\begin{figure}[t]
\centerline{\resizebox{5in}{4in}{
\includegraphics*{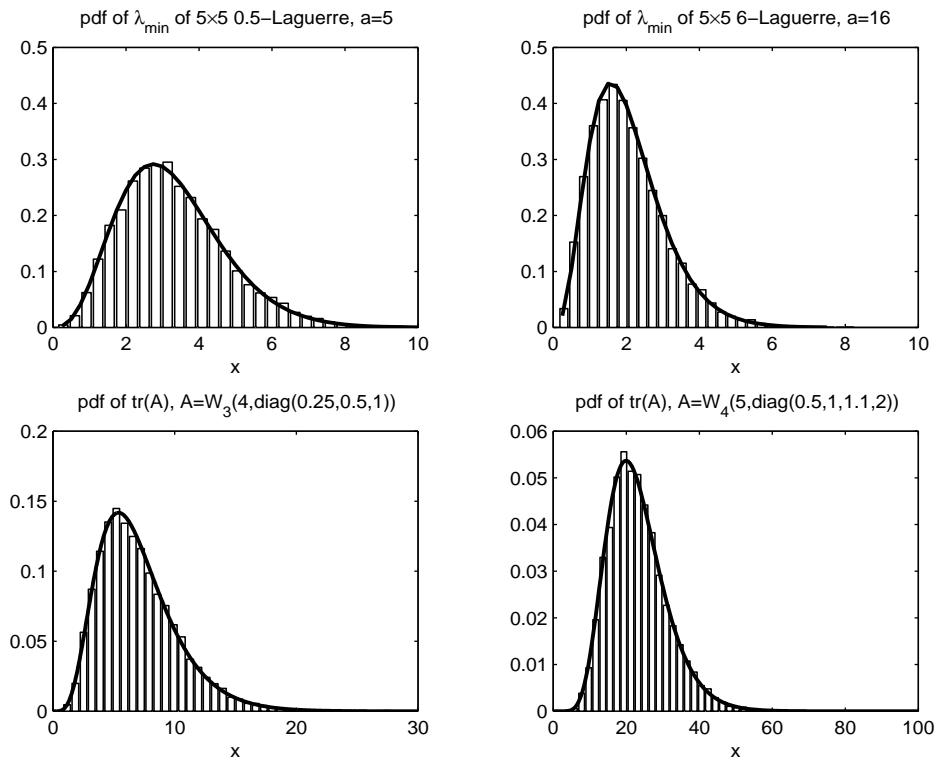}
}}
\caption{
The p.d.f.\ of $\lambda_{\min}$ of the $\beta$-Laguerre matrix
(top plots) and the p.d.f.\ of the trace of the Wishart matrix (bottom
plots).}
\label{fig_example2}
\end{figure}

Next we consider the distribution of the smallest eigenvalue of $L$ and 
that of $\tr\, A$.

If $c=a-\frac{\beta}{2}(n-1)-1$ is a nonnegative integer,
then the probability density function of the smallest eigenvalue of $L$
is (see, e.g.,~\cite[Thm.\ 10.1.1, p.\ 146]{dumitriuthesis}):
\begin{equation}
   f(x)=x^{cn}\cdot e^{-\frac{nx}{2}}\cdot
      \phantom{^\beta}_2F_0^{\mbox{\tiny$(2/\beta)$}}
\big(-c,\beta \textstyle{\frac{n}{2}+1;-\frac{2}{x}I_{n-1}}\big).
\label{example1}
\end{equation}
Since $-c$ is a nonpositive integer, the series expansion of 
$\!\!\phantom{^{\mbox{\tiny$\beta$}}}_2F_0^{\mbox{\tiny$(2/\beta)$}}$ 
in (\ref{example1}) terminates, even though it diverges in general.

The probability density function of $\tr\, A$
is \begin{eqnarray}
f(u)&=&\det(\lambda^{-1}\Sigma)^{-l/2}\sum_{k=0}^\infty
\frac{g_{\frac{ln}{2}+k,2\lambda}(u)}{k!}
\sum_{\kappa\vdash
k}\big({\textstyle\frac{l}{2}}\big)_\kappa^{\mbox{\tiny$(2)$}} \cdot
C_\kappa^{\mbox{\tiny$(2)$}}(I-\lambda\Sigma^{-1})
\nonumber
\\
&=&\det(\lambda^{-1}\Sigma)^{-l/2}\sum_{k=0}^\infty
 g_{\frac{ln}{2}+k,2\lambda} (u)
\sum_{\kappa\vdash k}
\big({\textstyle\frac{l}{2}}\big)_\kappa^{\mbox{\tiny$(2)$}} 2^k
j_\kappa^{-1}
\cdot
J_\kappa^{\mbox{\tiny$(2)$}}(I-\lambda\Sigma^{-1}),
\label{trA}
\end{eqnarray}
where
$$
g_{r,2\lambda}(u)=\frac{e^{-u/2\lambda}
u^{r-1}}{(2\lambda)^r\Gamma(r)}, \;\;\;\;\;(u>0),
$$
and $\lambda$ is arbitrary~\cite[p.~341]{muirhead}. We follow the
suggestion by Muirhead to use
$\lambda=2\lambda_1\lambda_l/(\lambda_1+\lambda_l)$, where $\lambda_1$ and
$\lambda_l$ are the largest and smallest eigenvalues of $\Sigma$, 
respectively.  Although
the expression~(\ref{trA}) is not a hypergeometric function of a matrix
argument, its truncation for $|\kappa|\le m$ has the form~(\ref{GX}),
and is computed analogously.

We plot~(\ref{example1}) and~(\ref{trA}) on Figure~\ref{fig_example2} and
compare the theoretical predictions of these formulas with experimental
data.

\subsection{Performance Results}
\label{sec_numexp_complexity}

\begin{figure}[t]
\centerline{\resizebox{5in}{2.5in}{
\includegraphics*{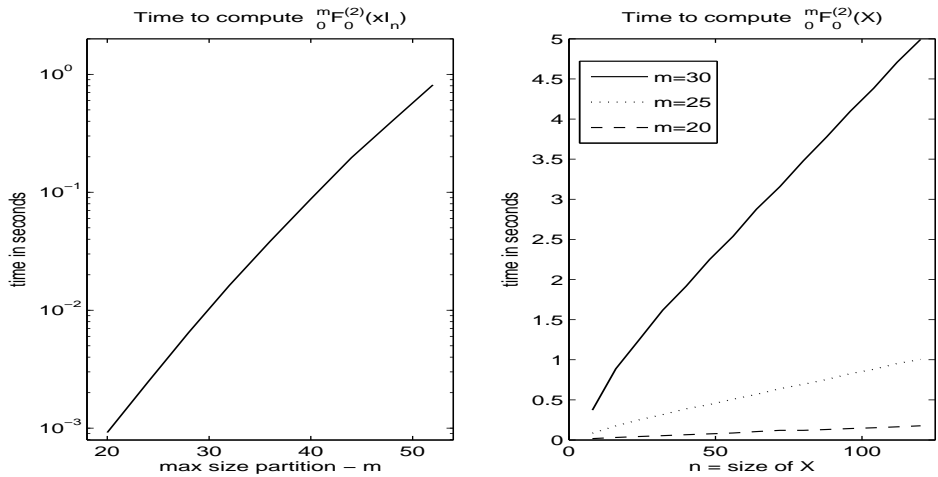}
}}
\caption{The performance of Algorithms~\ref{alg_hgi} (left) and
\ref{alg_hg} (right).}
\label{fig_example3}
\end{figure}

In Figure~\ref{fig_example3} we demonstrate the efficiency of
Algorithms \ref{alg_hgi} and \ref{alg_hg} on an 1.8GHz Intel Pentium 4
machine.

In the left plot we present the performance data for
Algorithm~\ref{alg_hgi} (whose complexity is independent of the size $n$ of
the matrix $X=xI$).  Its efficiency is evident---we need to sum beyond
partitions of size $m=52$ before this algorithm takes a full second (for
reference, the $k!$ in the denominator of the $\kappa$-term
in~(\ref{hg1}) then reaches up to $52!\approx 8\cdot 10^{67}$).

Algorithm~\ref{alg_hg} is also very efficient.  The right plot of
Figure~\ref{fig_example3} demonstrates clearly its linear complexity in
$n$.  It also takes at most a few seconds on matrices of size $n\le 120$
and partitions of size $m\le 30$.

\section{Conclusions and Open Problems}
\label{sec_concl}

We have presented new algorithms for computing the truncation of the 
hypergeometric function of a matrix argument. They exploit the 
combinatorial properties of the Pochhammer symbol and the Jack function to 
achieve remarkable efficiency, and have lead to 
new results~\cite{absiledelmankoev, edelmansutton04}.

Several problems remain open, among them 
automatic detection of convergence.
The $\kappa$-term in~(\ref{hg3}) does approach zero as
$|\kappa|\rightarrow\infty$, but it need not monotonically decrease.
Although we have
$$
\prod_{j=1}^{\kappa_i-1}
\frac
{(g_j-\alpha)e_j} {g_j(e_j+\alpha)}
\cdot
\prod_{j=1}^{i-1}
\frac {l_j-f_j} {l_j+h_j}
<1
$$
in~(\ref{Qkappaupdate}), it is not always true that
$Q_\kappa J_\kappa^{\al}\le Q_{\kappa_{(i)}}
J_{\kappa_{(i)}}^{\al}$, and it is unclear how to tell when
convergence sets in.

Another open problem is to determine the best way to truncate the
series~(\ref{hg1}).  Our choice to truncate it for $|\kappa|\le m$ seems to
work well in practice, but one can imagine selecting a partition $\lambda$
and truncating for $\kappa\le\lambda$ instead of, or in addition to
$|\kappa|\le m$.

\section*{Acknowledgements}
We thank Brian Sutton and Per-Olof Persson for
several discussions that resulted in the simplification of the
implementation of our algorithms, as well as the anonymous referee for the 
useful comments, which lead to improvements in the exposition.

\bibliographystyle{amsplain}

\providecommand{\bysame}{\leavevmode\hbox to3em{\hrulefill}\thinspace}
\providecommand{\MR}{\relax\ifhmode\unskip\space\fi MR }
\providecommand{\MRhref}[2]{%
  \href{http://www.ams.org/mathscinet-getitem?mr=#1}{#2}
}
\providecommand{\href}[2]{#2}

\end{document}